\documentclass[letterpaper]{article}
\usepackage{amsmath, hyperref}
\usepackage{graphicx}
\usepackage{amsfonts}
\usepackage{amssymb}

\newcommand{\comments}[1]{}

\newtheorem{theorem}{Theorem}

\newtheorem{assumption}[theorem]{Assumption}

\newtheorem{definition}[theorem]{Definition}

\newtheorem{lemma}[theorem]{Lemma}

\newtheorem{remark}[theorem]{Remark}

\title{\bf Constructions of Strict Lyapunov Functions for\\
Discrete Time and Hybrid Time-Varying Systems$\,
$\footnote{Corresponding Author: Fr\'ed\'eric Mazenc.}}
\author{Michael Malisoff$\,
$\footnote{Department of Mathematics; Louisiana State University;
Baton Rouge, LA 70803-4918; malisoff@lsu.edu.  Supported by NSF
Grant 0424011.  The author thanks Chaohong Cai, Rafal Goebel,
Ricardo Sanfelice, and Andrew Teel for illuminating discussions at
the 2005 American Control Conference in Portland, Oregon.}\and
Fr\'ed\'eric Mazenc$\, $\footnote{Projet MERE INRIA-INRA; UMR
Analyse des Syst\`emes et Biom\'etrie INRA; 2 pl. Viala; 34060
Montpellier, France; Frederic.Mazenc@ensam.inra.fr.}}
\begin{document}
\maketitle

\begin{abstract} We provide explicit
closed form expressions for strict Lyapunov functions for
time-varying discrete time  systems.  Our Lyapunov functions are
expressed in terms of known nonstrict Lyapunov functions for the
dynamics and finite sums of persistency of excitation parameters.
This provides a discrete time analog of our previous continuous time
Lyapunov function constructions. We also construct explicit strict
Lyapunov functions for systems satisfying nonstrict discrete time
analogs of the conditions from Matrosov's Theorem. We use our
methods to build strict Lyapunov functions for time-varying hybrid
systems that contain mixtures of
continuous and discrete time evolutions.\\[1em]
\end{abstract}

\section{Introduction}
\label{sec:intro} The theory of Lyapunov functions plays a
fundamental role in modern nonlinear robustness analysis and
controller design \cite{AlS99, AS99, ASW00, KSW01, MGSW05a, MGSW05b,
 MM05,
M03, MN06}. In many applications, it is essential to have explicit
closed form expressions for a strict Lyapunov function. This is
especially the case when one wishes to design stabilizing feedbacks,
which are often expressed in terms of the Lie derivatives of
Lyapunov functions in the directions of the vector fields that
define the system evolution. The classical converse Lyapunov
function theorem asserts that systems that are stable in an
appropriate sense also admit strict Lyapunov functions \cite{BR01}.
However, the Lyapunov functions provided by the theory are not
closed form explicit expressions since they involve infinite sums or
improper integrals or optimal control value functions and so do not
lend themselves to applications. Moreover, whereas most of the known
explicit Lyapunov function constructions are for time-invariant
systems, it is well appreciated that time-invariant systems are
often inadequate for engineering practice.  For example, there are
many applications where the dynamics cannot be stabilized by
time-invariant feedback but can be stabilized using time-varying
controllers \cite{C92, MSP99, Sa90, So90}. Time-varying systems are
also ubiquitous in tracking. While some methods for building
Lyapunov functions for time-varying systems are known, general
methods for  constructing {\em explicit closed form} Lyapunov
functions for time-varying discrete and hybrid systems are not
available. Hybrid systems are ubiquitous in science and engineering
\cite{VS00}. Hence, the construction of explicit Lyapunov functions
for time-varying systems presents significant challenges that are of
considerable ongoing research interest.

One  recently developed and  powerful approach to  this problem
involves constructing strict Lyapunov functions in terms of given
{\em non}strict Lyapunov functions for the system; see for instance
\cite{A99, FP00, M03, MN06}. By a nonstrict Lyapunov function, we
roughly mean a function that is positive definite and radially
unbounded and that has a negative semi-definite derivative along all
solutions; see Section \ref{sec:two} for precise definitions. The
advantage of this strictification approach is that in many
applications, a nonstrict Lyapunov function is readily available
through backstepping or physical considerations \cite{MM05}. For
continuous time systems, the strictification approach has been
applied to rotating rigid bodies, robot manipulators, and other
important engineering applications. This suggests the possibility of
extending the strictification approach by constructing explicit
closed form Lyapunov functions for discrete time time-varying
systems or, more generally, hybrid time-varying systems containing
both continuous and discrete time evolutions. The purpose of this
work is to  show that both of the extensions are indeed possible.

In Section \ref{sec:two}, we provide the relevant definitions of
strict Lyapunov functions and the necessary formalism of hybrid
systems, hybrid time domains, and hybrid trajectories.  In Section
\ref{sec:3}, we show how to construct explicit closed form strict
Lyapunov functions for time-varying discrete time nonlinear systems
in terms of given nonstrict Lyapunov functions. This provides a
discrete time analog of \cite{M03}  as well as a more explicit
Lyapunov function construction than the known discrete time
constructions that involve infinite sums of persistency of
excitation (PE) parameters \cite{NL04}.  We also build Lyapunov
functions for time-varying systems under appropriate versions of the
assumptions from Matrosov's Theorem thus providing a discrete time
analog of the results  \cite{MN06}  on continuous time systems
satisfying the Matrosov conditions.  In Section \ref{sec:4}, we
merge our results with the known continuous time analogs to
construct explicit closed form Lyapunov functions for time-varying
hybrid systems, under appropriate hybrid analogs of the PE or
Matrosov conditions. To our knowledge, this provides the first
general method for explicitly constructing Lyapunov functions for
 general time-varying nonlinear hybrid systems. In Section
\ref{sec:5}, we prove our theorems.  We provide some examples
covered by our results in Section \ref{sec:examples}, and we close
in Section \ref{sec:6} with some remarks about possible extensions.

\section{Definitions, Assumptions, and Lemmas} \label{sec:two}

We let ${\cal K}_\infty$ denote the set of all continuous functions
$\rho:[0,\infty)\to[0,\infty)$ for which (i) $\rho(0)=0$ and  (ii)
$\rho$ is strictly increasing and unbounded. Note that ${\cal
K}_\infty$ is closed under inverse and composition; i.e., if
$\rho_1,\rho_2\in {\cal K}_\infty$, then  $\rho^{-1}_1,
\rho_1\circ\rho_2\in {\cal K}_\infty$. We let ${\cal KL}$ denote the
class of all continuous functions $\beta:[0,\infty)\times
[0,\infty)\to[0,\infty)$ for which (I)  $\beta(\cdot, t)\in {\cal
K}_\infty$ for each $t\ge 0$, (II) $\beta(s,\cdot)$ is
non-increasing for each $s\ge 0$, and (III) $\beta(s,t)\to 0$ as
$t\to +\infty$ for each $s\ge 0$.  We let $\cal{KLL}$  denote the
set of all functions $\beta: [0,\infty)\times [0,\infty)\times
[0,\infty)\to [0,\infty)$ such that for each $\bar t\ge 0$, the
functions $(s,t)\mapsto\beta(s,t,\bar t)$ and $(s,t)\mapsto
\beta(s,\bar t, t)$ are of class $\mathcal{KL}$. When we say that a
function $\rho$ is {\em smooth} (a.k.a. $C^1$), we mean it is
continuously differentiable, written $\rho\in C^1$. (For functions
$\rho$ defined on $[0,\infty)$, we interpret $\rho'(0)$ as a
one-sided derivative, and continuity of $\rho'$ at $0$ as one-sided
continuity.)

We set ${\mathbb Z}_{\ge 0}=\{0,1,2,\ldots\}$, we let ${\mathbb
R}^n$ denote the set of all real $n$-tuples, and we use $|\cdot|$ to
denote the usual Euclidean norm.
   We say that a
function $\Theta:{\mathbb R}^n\times [0,\infty)\times {\mathbb
Z}_{\ge 0}\to {\mathbb R}: (x,t,k)\mapsto \Theta(x,t,k)$ (which may
be independent of $t$ or $k$) is {\em uniformly state-bounded} and
write $\Theta\in \mathcal{USB}$ provided there exists $\mu\in
{\mathcal K}_\infty$ such that $|\Theta(x,t,k)|\le \mu(|x|)$ for all
$x\in {\mathbb R}^n$, $t \ge 0$, and $k\in {\mathbb Z}_{\ge 0}$.
More generally, a vector valued function $H:{\mathbb R}^n\times
[0,\infty)\times {\mathbb Z}_{\ge 0}\to {\mathbb R}^n:
(x,t,k)\mapsto H(x,t,k)$ is of {\em class $\mathcal{USB}$}, written
$H\in \mathcal{USB}$,
 provided
$(x,t,k)\mapsto |H(x,t,k)|$ is of class $\mathcal{USB}$.
 Following \cite{MM05}, we also say $\Theta$  is {\em uniformly
proper and positive definite (UPPD)} and write $\Theta\in
\mathcal{UPPD}$ provided there are $\alpha_1,\alpha_2\in {\mathcal
K}_\infty$ such that $\alpha_1(|x|) \le  \Theta(x,t,k) \le
\alpha_2(|x|)$ for all $x\in {\mathbb R}^n$, $t\in [0,\infty)$,  and
$k\in {\mathbb Z}_{\ge 0}$. We say $\Theta$ is
$(\omega_1,\omega_2)$-{\em periodic} provided $\omega_1\in
[0,\infty)$ and $\omega_2\in {\mathbb Z}_{\ge 0}$ satisfy
\[\Theta(x,t-\omega_1, k-\omega_2)=\Theta(x,t,k)\; \; \; \; \forall (x,t,k)\in
{\mathbb R}^n\times [0,\infty)\times {\mathbb Z}_{\ge 0}.\]  When
$\Theta$ is independent of $t$ (resp., $k$), we define
$\omega_2$-periodicity (resp., $\omega_1$-periodicity) analogously.
A  continuous function defined on a subset of Euclidean space that
includes $0$ and valued in $(-\infty, 0]$ is {\em negative
semi-definite} provided it is  zero at zero. A continuous function
$\alpha$ defined on a subset of Euclidean space and valued in
$[0,\infty)$ is {\em positive definite} provided $\alpha$ is zero
only at zero in which case we also write $\alpha\in \mathcal{PD}$.

We study the stability properties of the discrete time fully
nonlinear time-varying system \begin{equation} \label{dis1}
x_{k+1}=F(x_k,k)
\end{equation}
where we always  assume $F\in \mathcal{USB}$.  We also study
continuous time time-varying systems
\begin{equation} \label{con1} \dot x=G(x,t)
\end{equation}
where $G\in \mathcal{USB}$ is locally Lipschitz.  We always assume
(\ref{con1}) is {\em forward complete}, meaning for each $x_o\in
{\mathbb R}^n$ and $t_o\ge 0$, there is a unique solution $t\mapsto
\phi(t,t_o, x_o)$ for (\ref{con1}) defined on $[t_o,\infty)$ that
satisfies $\phi(t_o, t_o,x_o)=x_o$. We interpret the solutions of
(\ref{con1}) in the generalized Lebesgue almost all (a.a.) sense. We
also use $k\mapsto \phi(k,k_o, x_o)$ to denote the discrete time
solution of (\ref{dis1}) satisfying $\phi(k_o, k_o, x_o)=x_o$
whenever this would not lead to confusion.  Given a function
$V:{\mathbb R}^n\times [0,\infty)\times {\mathbb Z}_{\ge 0}\to
{\mathbb R}: (x,t,k)\mapsto V(x,t,k)$, we set
\[
\renewcommand{\arraystretch}{1.25}
\Delta_kV(x,t,k):=V(F(x,k),t,k+1)-V(x,t,k),\; \; \; \;  {\mathcal
D}V(x,t,k):=\frac{\partial V}{\partial t}(x,t,k)+\frac{\partial
V}{\partial x}(x,t,k)G(x,t)
\]
assuming $(x,t)\mapsto V(x,t,k)$ is also smooth for each $k\in
\mathbb{Z}_{\ge 0}$ in the definition of ${\mathcal D}V$.  In our
analysis of (\ref{dis1}), $V$ will generally not depend on $t$ but
we need to allow its dependence on $t$ in our discussion of hybrid
systems.

\begin{definition}\label{disdef}(a)  Let $V:{\mathbb R}^n\times
{\mathbb Z}_{\ge 0}\to {\mathbb R}$ be of class $\mathcal{UPPD}$.
We call $V$ a {\em (strict) Lyapunov function} for (\ref{dis1})
provided there exists $\alpha_3\in \mathcal{PD}$ such that:
\begin{equation}
\label{olv} \Delta_k V(x,k) \; \le\;  -\alpha_3(|x|)\; \; \; \; \;
\; \forall x\in {\mathbb R}^n \,  \; \&\, \;  k\in {\mathbb Z}_{\ge
0}.
\end{equation}
(b) We say that (\ref{dis1}) is {\em globally asymptotically stable
(GAS)} provided there exists $\beta\in \mathcal{KL}$ such that for
all $x_o\in {\mathbb R}^n$ and $k_o\in {\mathbb Z}_{\ge 0}$, we have
$|\phi(k,k_o, x_o)|\le \beta(|x_o|,k-k_o)$ for all $k\ge k_o$.
\end{definition}
The corresponding Lyapunov function and GAS definitions for
(\ref{con1}) are obtained from Definition \ref{disdef} by replacing
$k$ with $t$, $\mathbb{Z}_{\scriptscriptstyle \ge 0}$ with
$[0,\infty)$,
 and $\Delta_kV$ with
${\mathcal D}V$. Notice that we do not require $\alpha_3$ to be of
class ${\mathcal K}_\infty$. For the special case where $\alpha_3\in
{\mathcal K}_\infty$, the existence of a discrete time Lyapunov
function $V$ is known to imply that (\ref{dis1}) is GAS since then
$\Delta_kV(x,k)\le -\alpha(V(x,k))$ everywhere with
$\alpha:=\alpha_3\circ\alpha^{-1}_2\in \mathcal{K}_\infty$ and
$\alpha_2$ as in the $\mathcal{UPPD}$ condition on $V$ \cite[Theorem
8]{NTS99}.
 Furthermore, by replacing the function $\alpha_3(|x|)$ in
(\ref{olv}) by the smaller function $\Theta(V(x,k)) :=
\displaystyle\min\{\alpha_3(s):\alpha_2^{-1}(V(x,k))\le s \le
\alpha_1^{-1}(V(x,k))\}$, we are in a situation where Lemma \ref{tl}
below applies and straightforwardly implies that one can construct a
Lyapunov function satisfying (\ref{olv}) with a new function
$\alpha_3$ of class ${\mathcal K}_\infty$. Combining this fact with
the stability result from \cite{NTS99}, we get:
\begin{lemma}
\label{impl} If (\ref{dis1}) admits a strict Lyapunov function, then
it is GAS.
\end{lemma}
We also use the following persistency of excitation (PE) notions
from \cite{MM05, NL04}:

\begin{definition}\label{pedef}
(a) We say that a bounded function $p:{\mathbb Z}_{\ge 0}\to
[0,\infty)$ is of {\em discrete PE type with parameters $l$ and
$\delta$} and write $p\in \mathcal{P}_{\rm dis}(l,\delta)$ provided
$l\in {\mathbb Z}_{\ge 0}$ and $\delta>0$ are such that
\begin{equation}
\label{ped} \displaystyle\sum_{i = k - l}^{k} p(i) \; \geq\;
\delta\; \; \; \; \forall k\in {\mathbb Z}_{\ge 0}.
\end{equation}
(b) We say that a bounded continuous function $q:[0,\infty)\to
[0,\infty)$ is of {\em continuous PE type with parameters  $\tau$
and $\varepsilon$} and write $q\in \mathcal{P}_{\rm
cts}(\tau,\varepsilon)$ provided $\tau\ge 0$ and $\varepsilon>0$ are
such that
\begin{equation}
\label{pec} \int_{t-\tau}^tq(r)\, dr\; \ge\;  \varepsilon\;  \; \;
\forall t\ge 0.\end{equation} (c) We set $\mathcal{P}_{\rm
dis}=\bigcup\{ \mathcal{P}_{\rm dis}(l,\delta): l\in {\mathbb
Z}_{\ge 0}, \delta>0\}$ and  $\mathcal{P}_{\rm cts}=\bigcup\{
\mathcal{P}_{\rm cts}(\tau,\varepsilon): \tau\ge 0,
\varepsilon>0\}$.
\end{definition}

Elements of $\mathcal{P}_{\rm dis}$ and $\mathcal{P}_{\rm cts}$ are
called {\em PE parameters} and arise in a variety of contexts, e.g.,
$q(t)=\sin^2(t)$ as well as cases where $q$ can be null on intervals
of arbitrarily large length \cite{MM05}. The following lemma follows
from a simple change of variables, a Fubini Theorem argument (as was
used in \cite{MM05}), and the formula $1+2+\ldots+m=m(m+1)/2$:

\begin{lemma}\label{fubl} Let $l\in \mathbb{Z}_{\ge 0}$, let
$\tau,\varepsilon, \delta>0$, and let $p\in {\mathcal P}_{\rm
dis}(l,\delta)$ and $q\in \mathcal {P}_{\rm cts}(\tau,\varepsilon)$
be bounded above by $\bar p$ and $\bar q$ respectively.  Define the
functions $S:{\mathbb Z}_{\ge 0}\to [0,\infty)$ and
$R:[0,\infty)\to[0,\infty)$ by
\begin{equation}\label{skdef}
S(k):=\displaystyle\sum_{s = k - l}^{k} \displaystyle\sum_{j =
s}^{k}\, p(j),\; \; \; \; R(t):=\int_{t-\tau}^t\int_z^tq(\nu)\,
d\nu\, dz.
\end{equation}
Then $S(k)\le \bar p(l+1)^2$ and $R(t)\le \tau^2\bar q/2$ hold for
all $k\in \mathbb{Z}_{\ge 0}$ and $t\ge 0$.  If $p$ is $l$-periodic,
then so is $S$.  If $q$ is $\tau$-periodic, then so is $R$.
\end{lemma}

We next recall the hybrid system tools developed in  \cite{CTG05,
C04}, generalized to time-varying systems. For simplicity, we only
consider singleton valued dynamics although our results carry
through in the more general setting of difference and differential
inclusions. Given sets $C,D\subseteq {\mathbb R}^n$ and $F$ and $G$
satisfying the assumptions above, the corresponding {\em hybrid
dynamical system} is defined to be the formal  object
\begin{equation}
\label{hybsy} {\mathcal H}:= \left\{
\begin{array}{lcll}
\dot x&=&G(x,t), & x\in C\\
x_{k+1}&=&F(x_k,k),&  x_k\in D
\end{array}.
\right.\end{equation} A {\em compact hybrid time domain} is  a
subset $E\subset [0,\infty)\times \mathbb{Z}_{\ge 0}$ of the form
$\cup_{k=0}^{K-1}([t_k, t_{k+1}]\times\{k\})$ for some finite
sequence $0\le t_o\le t_1\le \ldots \le t_K$.  A {\em hybrid time
domain} is a set $E\subset [0,\infty)\times \mathbb{Z}_{\ge 0}$ with
the property that for all $(T,K)\in E$, the intersection
$E\cap([0,T]\times\{0,1,\ldots, K\})$ is a compact hybrid time
domain. A {\em hybrid arc} is a function $x(t,k)$ defined on a
hybrid time domain ${\rm dom}(x)$ such that $t\mapsto x(t,k)$ is
locally absolutely continuous for each $k$. A {\em hybrid
trajectory} of (\ref{hybsy}) is  a hybrid arc $x(t,k)$ that
satisfies the following:
\begin{itemize}
\item[$(S_1)$] For all $k\in {\mathbb Z}_{\ge 0}$ and a.a. $t$
such that $(t,k)\in {\rm dom}(x)$, we have $x(t,k)\in C$ and
$\frac{\partial}{\partial t} x(t,k)=G(x(t,k),t)$.
\item[$(S_2)$] For all $(t,k)\in {\rm dom}(x)$ such that  $(t,k+1)\in{\rm dom}(x)$, we have $x(t,k)\in D$ and
$x(t,k+1)=F(x(t,k),k)$.\end{itemize}

Notice that $E$ is a hybrid time domain provided it is a finite or
infinite union of sets of the form $[t_k,t_{k+1}]\times \{k\}$ with
$\{t_k\}$ nondecreasing in $[0,\infty)$, with a possible additional
`last' set having the form $[t_{k}, T)\times\{k\}$ with $T$ finite
or infinite. To keep our notation simple, we use $\cup_{k\in
J}([t_k, t_{k+1}]\times\{k\})$ to denote a generic hybrid time
domain with the understanding that (i) either $J=\mathbb{Z}_{\ge 0}$
or $J$ is a finite set of the form $\{0,1,2,\ldots, j_{\rm max}\}$
and (ii) $[t_k, t_{k+1}]$ may mean $[t_k, t_{k+1})$ if  $J$ is
finite and $k=j_{\rm max}$.  Notice that continuous time solutions
of (\ref{con1}) in $C$ and discrete time solutions of (\ref{dis1})
in $D$  starting with  $k=0$ correspond to hybrid trajectories of
(\ref{hybsy}) that have no switchings between the discrete and
continuous evolutions.

\begin{definition}\label{hyly}
(a) Let $V\in \mathcal{UPPD}$ be $C^1$ in $x$ and $t$.  We call $V$
a {\em (strict) Lyapunov function} for $\mathcal H$ provided there
exists $\alpha_3\in \mathcal{PD}$ such that the following hold for
all $t\ge 0$ and $k\in {\mathbb Z}_{\ge 0}$:
\begin{equation}
\label{hyly2} \Delta_kV(x,t,k)\; \le\; -\alpha_3(|x|)\; \; \;
\forall x\in D;\; \; \;  \; {\mathcal D}V(x,t,k)\; \le\;
-\alpha_3(|x|)\; \; \; \forall x\in C.
\end{equation}
If, in addition, there is a constant $r>0$ such that
\begin{equation}
\label{hyl2a} V(F(x,k), t, k+1)\; \le \; e^{-r}V(x,t,k)\; \; \;
\forall x\in D;\; \; \; \; {\mathcal D}V(x,t,k)\; \le\; -rV(x,t,k)\;
\; \; \forall x\in C,
\end{equation}
then we call $V$ an {\em exponential decay Lyapunov function} for
${\mathcal H}$. (b) We call $\mathcal H$ {\em globally
asymptotically stable (GAS)} provided there exists $\beta\in
\mathcal{KLL}$ such that: For each trajectory $x(t,k)$ of $\mathcal
H$ defined on any hybrid time domain  $\cup_{k\in J}([t_k,
t_{k+1}]\times \{k\})$,  we have $|x(t,k)|\le \beta(|x(t_o,0)|, k,
t-t_k)$ for all $k\in J$ and all $t\in [t_k, t_{k+1}]$.
\end{definition}

\begin{lemma}\label{implh} If $\mathcal H$ admits a Lyapunov
function, then it is GAS.\end{lemma}

To prove this lemma, first note that since $\alpha_3$ in
(\ref{hyly2}) is independent of $k$,  standard arguments (e.g. those
in \cite{S89} applied with $a(x):=\alpha_3(|x|)$) provide
$\beta_1\in \mathcal{KL}$ such that for each hybrid trajectory
$x(t,j)$ defined on a hybrid time domain $\cup_{k\in J}([t_k,
t_{k+1}]\times\{k\})$ and satisfying any initial condition
$x(t_o,0)=x_o$, we have
\begin{equation}
\label{beta1} |x(t,k)|\; \le\;  \beta_1(|x(t_k,k)|,t-t_k)\; \; \; \;
\forall k\in J,\; \; t\in [t_k, t_{k+1}].\end{equation} Similarly,
since $\alpha_3$ in (\ref{hyly2}) is independent of $t$, and since
we can assume as above that $\alpha_3\in \mathcal{K}_\infty$, the
argument from \cite[Theorem 8]{NTS99}  provides $\beta_2\in
\mathcal{KL}$ such that
\begin{equation}
\label{beta2} |x(t_k,k)|\; \le\;  \beta_2(|x_o|, k)\; \; \; \;
\forall k\in J.\end{equation} In fact, $\beta_2$ can be constructed
using the decay conditions from (\ref{hyly2}) as follows. First note
that by arguing as in the proof of Lemma \ref{impl} above and
replacing $V$ with $\kappa\circ V$ for a suitable function
$\kappa\in \mathcal{K}_\infty$ in the discrete decay condition
without relabeling, we can find $\gamma\in \mathcal{K}_\infty$ such
that $\Delta_kV(x,t,k)\le -\gamma(V(x,t,k))$ when $x\in D$; see
Lemma \ref{tl} for the construction of $\kappa$.   Since $t\mapsto
V(x(t,k),t,k)$ decays on $(t_k,t_{k+1})$ for each $k$, we get
$V(x(t_{k+2}, k+1),t_{k+2},k+1)\le  V(x(t_{k+1},
k+1),t_{k+1},k+1)=V(F(x(t_{k+1}, k),k),t_{k+1},k+1)$ so
\[\begin{array}{rcl}
V(x(t_{k+2}, k+1),t_{k+2},k+1)-V(x(t_{k+1}, k),t_{k+1},k)&\le &
\Delta_kV(x(t_{k+1}, k),t_{k+1},k)\\& \le & -\gamma(V(x(t_{k+1},
k),t_{k+1},k))
\end{array}\]
everywhere.  Applying \cite[Theorem 8]{NTS99} to the function
$k\mapsto V(x(t_{k+1}, k),t_{k+1},k)$ and recalling that $V$ is
uniformly proper and positive definite
 gives $\tilde\beta_2\in
\mathcal{KL}$ (not depending on the choice of the trajectory) such
that $|x(t_{k+1},k)|\le \tilde \beta_2(|x(t_1,0)|,k)$ for all $k$.
Choosing $\alpha_1,\alpha_2\in \mathcal{K}_\infty$ such that
$\alpha_1(|x|)\le V(x,t,k)\le \alpha_2(|x|)$ everywhere, the
discrete time decay condition in (\ref{hyly2}) gives
\[|x(t_{k+1},k)|\ge \alpha^{-1}_2\circ V(x(t_{k+1}, k),t_{k+1},k)\ge
\alpha^{-1}_2\circ V(x(t_{k+1}, k+1),t_{k+1},k+1)\ge
\alpha^{-1}_2\circ\alpha_1(|x(t_{k+1}, k+1)|)\] for all $k\in J$.
Similarly, the continuous time decay condition in (\ref{hyly2})
gives \[|x(t_1,0)|\le \alpha^{-1}_1\circ V(x(t_1,0),t_1,0)\le
\alpha^{-1}_1\circ V(x(t_o,0),t_o,0)\le \alpha^{-1}_1\circ
\alpha_2(|x_o|).\] We can therefore satisfy (\ref{beta2}) by taking
$\beta_2(s,k):=\alpha^{-1}_1\circ\alpha_2\circ\tilde
\beta_2(\alpha^{-1}_1\circ\alpha_2(s),k)+s/(k+1)$, where the
additional term $s/(k+1)$ is used to account for the case $k=0$.
Combining (\ref{beta1})-(\ref{beta2}) shows we can satisfy the
requirements of
 Lemma
\ref{implh} using $\beta(s,t,k)=\beta_1(\beta_2(s,k),t)$.

\section{Statement of Results on Discrete-Time Systems}
\label{sec:3} \subsection{Strictifying Persistence of Excitation
(PE) Decay Estimates} We begin by constructing explicit closed form
Lyapunov functions for discrete time systems in terms of nonstrict
Lyapunov functions and appropriate PE parameters $p\in {\mathcal
P}_{\rm dis}$.  For an alternative construction, involving infinite
sums of PE parameter values, see \cite{NL04}. We prove the following
in
 Section \ref{sec:5}:
\begin{theorem}
\label{disthm} Let $l\in {\mathbb Z}_{\ge 0}$, $\delta>0$,  $p\in
{\mathcal P}_{\rm dis}(l,\delta)$, $V\in \mathcal{UPPD}$, and
 $\Theta \in \mathcal{PD}$ satisfy
\begin{equation}
\label{PEdecay} \Delta_kV(x,k)\; \le\; -p(k+1)\Theta(V(x,k))\; \; \;
\forall x\in {\mathbb R}^n\, \; \& \, \; k\in {\mathbb Z}_{\ge 0}.
\end{equation}
Then one can construct $\kappa,\gamma\in {\mathcal K}_\infty$ such
that
\begin{equation}
\label{uformula} U(x,k)\; :=\;
\kappa(V(x,k))+\frac{\gamma(V(x,k))}{4(l+1)}\displaystyle\sum_{s = k
- l}^{k} \displaystyle\sum_{j = s}^{k} p(j)
\end{equation}
is a strict Lyapunov function for (\ref{dis1}), so (\ref{dis1}) is
GAS. If  $p$ and $V$ are also both $l$-periodic in $k$, then so is
$U$.

\end{theorem}

\begin{remark}
A key feature in (\ref{PEdecay}) is that the PE condition on $p$
allows $p(k+1)=0$ for some values of $k$ in which case we could have
$\Delta_kV(x,k)=0$.
 An additional novel feature of
Theorem \ref{disthm} is that we do not require the gain function
$\Theta$ in (\ref{PEdecay}) to be of class ${\mathcal K}_\infty$.
This properness of the gain function was required in \cite{MM05,
NL04}. Our proof of Theorem \ref{disthm} will show that we can take
$\kappa(s)\equiv s$ if $\Theta\in \mathcal{K}_\infty$.
\end{remark}

\subsection{Lyapunov Function Constructions Under Matrosov
Conditions}

Recall the definitions of $\mathcal{UPPD}$ and $\mathcal{USB}$ from
Section \ref{sec:two}. We explicitly construct a  Lyapunov function
for discrete time systems (\ref{dis1}) satisfying the following
analog of the Matrosov Theorem conditions from \cite{MN06}:

\begin{assumption}\label{as1}  There exist
$V_1:{\mathbb R}^n\times {\mathbb Z}_{\ge 0}\to [0,\infty)$ of class
$\mathcal{UPPD}$, $V_2:{\mathbb R}^n\times {\mathbb Z}_{\ge 0}\to
{\mathbb R}$ of class $\mathcal{USB}$, a function $\phi_2\in
{\mathcal K}_\infty$, nonnegative functions $N_1, N_2\in
\mathcal{USB}$, a function $\chi:{\mathbb R}^n\times
[0,\infty)\times {\mathbb Z}_{\ge 0}\to {\mathbb R}$, a positive
increasing function $\phi_1$, a positive definite function $W$, and
$p\in {\cal P}_{\rm dis}$ such that
\[\renewcommand{\arraystretch}{1.25}
\begin{array}{l}
 \Delta_k V_{1}(x,k) \leq - N_1(x,k),
\\[.5em]
 \Delta_k
V_{2}(x,k) \leq - N_2(x,k) + \chi(x,N_1(x,k),k),\\[.5em]
 |\chi(x,N_1(x,k),k)| \leq
\phi_1(|x|)\phi_2(N_1(x,k)),\; \; {\rm and}\; \; N_1(x,k) + N_2(x,k)
\geq p(k+1) W(x)
\end{array}\]
hold  for all $x\in {\mathbb R}^n$ and $k\in {\mathbb Z}_{\ge 0}$.
\end{assumption}   Notice
that we allow $V_2$ to take both positive and negative values. In
Section \ref{sec:5}, we prove:

\begin{theorem}
\label{matrthm} If (\ref{dis1}) satisfies Assumption \ref{as1}, then
one can construct an explicit strict Lyapunov function for
(\ref{dis1}).  In particular, (\ref{dis1}) is GAS.\end{theorem}

\section{Statement of Results on Hybrid Systems}
\label{sec:4}
\subsection{Hybrid Persistency of Excitation Estimates}
We next extend Theorem \ref{disthm} to hybrid systems. To keep the
exposition simple, we  assume the gain functions $\Theta$ in
(\ref{PEdecay}) and its continuous analog are $\Theta(s)=s$, but the
extension to general positive definite $\Theta$ can be done using
similar arguments.  We prove the following in Section \ref{sec:5}:

\begin{theorem}\label{hybthm}
 Let $V\in \mathcal{UPPD}$ be $C^1$ in $x$ and $t$. Consider the hybrid
system ${\cal H}$ in (\ref{hybsy}), and let $\delta,\varepsilon,
\tau>0$ and $l\in {\mathbb Z}_{\ge 0}$ be given.
 Assume there exist $r\in
\mathcal{P}_{\rm dis}(l,\delta)$ and $q\in \mathcal{P}_{\rm
cts}(\tau,\varepsilon)$
 such that
\begin{equation}
\label{hyply2}
 V(F(x,k), t, k+1)\; \le e^{-r(k+1)}V(x,t,k)\; \; \;
\forall x\in D;\; \; \; \; \; \; \;
 {\mathcal D}V(x,t,k)\; \le\; -q(t)V(x,t,k)\; \; \;
\forall x\in C
\end{equation}
 hold for all $t\ge 0$ and $k\in{\mathbb Z}_{\ge 0}$.
Then
\begin{equation}
\label{vs}
V^\sharp(x,t,k)=\left[2+\frac{1}{4(l+1)}\displaystyle\sum_{s = k -
l}^{k} \displaystyle\sum_{j = s}^{k}\left(1-e^{-r(j)}\right)+
\frac{1}{\tau}\int_{t-\tau}^t\int_z^tq(\nu)\, d\nu\, dz
\right]V(x,t,k)\end{equation} is an exponential decay Lyapunov
function for $\mathcal H$ which is therefore GAS.  If in addition
$V$ is $(\tau,l)$-periodic and $r$ and $q$ are $l$-periodic and
$\tau$-periodic respectively, then $V^\sharp$ is also
$(\tau,l)$-periodic.
\end{theorem}

\begin{remark}\label{cases}
The preceding theorem covers continuous dynamics (by taking
$D=\emptyset$ and $C=\mathbb{R}^n$ with the understanding that the
term involving the double sum in $V^\sharp$ is not present) and
discrete dynamics (by taking $C=\emptyset$ and $D=\mathbb{R}^n$ in
which case the term involving the double integral in $V^\sharp$ is
not present).  See \cite{CTG05} for an alternative, nonexplicit
construction of a Lyapunov function for time-invariant hybrid
systems.
\end{remark}

\subsection{Hybrid Systems Satisfying Matrosov Conditions}

We next extend Theorem \ref{matrthm} to hybrid systems that satisfy
the following analog of Assumption \ref{as1}.

\begin{assumption}\label{as2}  There exist $V_1\in \mathcal{UPPD}$ and
$V_2\in \mathcal{USB}$ that are $C^1$ in $(x,t)$,   nonnegative
$N_1, N_2\in {\rm \mathcal{USB}}$, a function $\chi:{\mathbb
R}^n\times [0,\infty)^2\times {\mathbb Z}_{\ge 0}\to {\mathbb R}$, a
positive increasing $\phi_1$, and  a positive definite function $W$,
$p\in {\cal P}_{\rm dis}$, $\phi_2\in {\mathcal K}_\infty$, and
$q\in {\mathcal P}_{\rm cts}$ such that
\begin{enumerate}\addtolength{\itemsep}{-0.35\baselineskip}
\item For all $x\in D$, we have
$\Delta_k V_{1}(x,t,k)  \leq  - N_1(x,t, k)$, $\Delta_kV_{2}(x,t,k)
 \leq  - N_2(x,t,k) + \chi(x,N_1(x,t,k),t,k)$, and  $N_1(x,t,k) +
N_2(x,t,k)  \geq  p(k+1) W(x)$.
\item For all $x\in C$, we have
${\mathcal D}V_1(x,t,k) \le  -N_1(x,t,k)$, ${\mathcal D}V_2(x,t,k)
\leq - N_2(x,t,k) + \chi(x,N_1(x,t,k),t,k)$, and  $N_1(x,t,k) +
N_2(x,t,k)  \geq  q(t) W(x)$.
\item For all $x\in \mathbb{R}^n$, we have
$|\chi(x,N_1(x,t,k),t, k)|  \leq \phi_1(|x|)\phi_2(N_1(x,t,k))$.
\end{enumerate}
 hold  for all $t\ge 0$  and $k\in {\mathbb Z}_{\ge 0}$.
\end{assumption}

 Assumption \ref{as2} simply means  the discrete and continuous
parts of $\mathcal{H}$ satisfy the appropriate discrete and
continuous Matrosov conditions.
 It
reduces to Assumption \ref{as1} for discrete systems when
$C=\emptyset$ and $D=\mathbb{R}^n$ in which case its condition 2.
holds vacuously. Notice that we again do not require $V_2$ to be
nonnegative. In Section \ref{sec:5}, we prove:
\begin{theorem}
\label{hmatrthm} If $\mathcal H$ satisfies Assumption \ref{as2},
then one can construct an explicit closed form strict Lyapunov
function for $\mathcal H$. In particular, $\mathcal H$ is
GAS.\end{theorem}

\section{Proofs of Theorems}
\label{sec:5}

\subsection{Results on Discrete Systems}
\subsubsection{Proof of Theorem \ref{disthm}}
By minorizing $\Theta\in \mathcal{PD}$ without relabeling as in
\cite{MM06}, we assume in the sequel that $\Theta\in C^1$ is
nondecreasing on $[0,1]$ and nonincreasing on $[1,\infty)$. The next
technical lemma allows us to assume that $\Theta\in
\mathcal{K}_\infty$ in (\ref{PEdecay}):

\begin{lemma} \label{tl}Let $\Theta\in \mathcal{PD}$ be as above and $p\in
{\mathcal P}_{\rm dis}$.  Define $\mu:[0, \infty)\to [1,\infty)$,
$\kappa$, and $\chi$ by \begin{equation}\label{kappadef} \kappa(r):=
2\int_0^r \mu(z)\, dz,\; \; \; \; \chi(r):=\Theta(2r)\mu(r),\; \; \;
\; {\rm and}\; \; \; \;
\mu(r)= \left\{\begin{array}{lcl} 1+4r^2,& \! \! 0\le r\le 1/2\\
\displaystyle\frac{4\Theta(1)r}{\Theta(2r)},& 1/2\le r<\infty
\end{array}\right..
\end{equation} Let $\nu\in \mathcal{UPPD}$ satisfy $\Delta_k\nu(x,k) \leq -
p(k+1) \Theta(\nu(x,k))$ for all $x\in {\mathbb R}^n$ and $k\in
\mathbb{Z}_{\ge 0}$. Then $\kappa\in \mathcal{K}_\infty\cap C^1$,
$\chi\in \mathcal{K}_\infty$, and $V := \kappa(\nu)\in
\mathcal{UPPD}$ satisfies
\begin{equation}
\label{u2} \Delta_kV(x,k)\leq - p(k+1) \gamma(V(x,k))\; \; \; \;
\forall x\in {\mathbb R}^n\; \&\; k\in {\mathbb Z}_{\ge 0},
\end{equation}
where $\gamma\in \mathcal{K}_\infty$ is defined by
$\gamma(s):=\chi(\kappa^{-1}(s)/2)$.
\end{lemma}
To prove Lemma \ref{tl}, fix $x\in \mathbb{R}^n$ and $k\in
\mathbb{Z}_{\ge 0}$ and apply the Fundamental Theorem of Calculus to
$s\mapsto {\mathcal F}(s):=
\kappa\left(s\nu(F(x,k),k+1)+(1-s)\nu(x,k)\right)$ to write
$\Delta_kV(x,k)={\mathcal F}(1)-{\mathcal
F}(0)=\int_{\scriptscriptstyle 0}^{\scriptscriptstyle 1}{\mathcal
F}'(s)\, ds$ and so also
\begin{eqnarray*}
\Delta_kV(x,k)&=&\left[\displaystyle\int_{0}^{1}
\kappa'\left(s\nu(F(x,k),k+1)
+ (1-s)\nu(x,k)\right) ds\right] \left[\nu(F(x,k),k+1) - \nu(x,k)\right]\\
&\le&- p(k+1)\left[\displaystyle\int_{0}^{1}
\kappa'(s\nu(F(x,k),k+1) + (1-s)\nu(x,k))
ds\right]\Theta(\nu(x,k))\\
& \leq & - p(k+1)\left[\displaystyle\int_{0}^{1}
\kappa'((1-s)\nu(x,k)) ds\right]\Theta(\nu(x,k))
\\
& \leq & - p(k+1)\left[\displaystyle\int_{0}^{1/2}
\kappa'\left(\frac{1}{2}\nu(x,k) \right) ds\right]\Theta(\nu(x,k))
\; \; = \; \; - p(k+1)\mu\left(\frac{1}{2}\nu(x,k)\right)
\Theta(\nu(x,k))
\end{eqnarray*}
where the first inequality holds because $\kappa$ is nondecreasing
and the other inequalities used the fact that $\kappa'$ is
nondecreasing. The lemma now follows from our choices of $\gamma$
and $\chi$.

We can therefore assume  that $V$ satisfies (\ref{u2}) with
$\gamma\in {\mathcal K}_\infty$, possibly by replacing $V$ with
$\kappa(V)$ for $\kappa\in \mathcal{K}_\infty$ defined in
(\ref{kappadef}). Defining $S(k)$ as in (\ref{skdef})
 and
defining $U$  by (\ref{uformula}) with $\kappa(s)\equiv s$ therefore
gives \begin{eqnarray} \label{key1}
 \! \Delta_kU(x,k) & \! \! \! = & \! \! V(F(x,k),k+1) + \frac{S(k+1)}{4(l+1)}
\gamma(V(F(x,k),k+1)) - V(x,k)
 - \frac{S(k)}{4(l+1)}\gamma(V(x,k))\nonumber\\
& = & \! \! \Delta_k V(x,k) + \frac{1}{4(l+1)}
S(k+1)\Delta_k(\gamma\circ V)(x,k)+
\frac{1}{4(l+1)}\gamma(V(x,k))\left[S(k+1) - S(k)\right]\nonumber\\
&\le & \! \! \Delta_k V(x,k) + \frac{1}{4(l+1)}\gamma(V(x,k))
\left[S(k+1) - S(k)\right],
\end{eqnarray}
where the last inequality holds because  $\gamma$ is increasing, so
$\Delta_k(\gamma\circ V)(x,k)\le 0$. Note that
\begin{equation}\label{key2}
\renewcommand{\arraystretch}{1.1}
\begin{array}{rcl}
 S(k+1) - S(k) & = &
\displaystyle\sum_{s = k+1 - l}^{k}\displaystyle\sum_{j = s}^{k}
p(j) + (l+1) p(k+1) - \! \displaystyle\sum_{s = k - l}^{k}
\displaystyle\sum_{j = s}^{k} p(j)
\\
& = & \displaystyle\sum_{s = k - l}^{k}\displaystyle\sum_{j = s}^{k}
p(j) - \displaystyle\sum_{j = k-l}^{k} p(j) + (l+1) p(k+1)
 - \displaystyle\sum_{s = k - l}^{k} \displaystyle\sum_{j =
s}^{k} p(j)
\\
& = & - \displaystyle\sum_{j = k-l}^{k} p(j) + (l+1) p(k+1)\, .
\end{array}\end{equation}
Substituting (\ref{key2}) into (\ref{key1}) gives
\begin{equation}
\label{10}
\renewcommand{\arraystretch}{1.1}
\begin{array}{rcl}
\Delta_k U(x,k) & \leq & \Delta_k V(x,k) +
\displaystyle\frac{1}{4(l+1)}\gamma(V(x,k))\left((l+1)
p(k+1)-\displaystyle\sum_{j = k-l}^{k} p(j)\right)
\\
& \leq & \Delta_k V(x,k) +
\displaystyle\frac{p(k+1)\gamma(V(x,k))}{4} -
\displaystyle\frac{\gamma(V(x,k))}{4(l+1)}\displaystyle\sum_{j =
k-l}^{k} p(j)\; \le \;
 - \frac{\delta}{4(l+1)} \gamma(V(x,k)),
\end{array}
\end{equation}
where the last inequality follows from the PE property of $p$ and
(\ref{u2}). This,  the fact that $V\in \mathcal{UPPD}$, and the
global boundedness of $S(k)$ from Lemma \ref{fubl} shows that $U$ is
a Lyapunov function for (\ref{dis1}).  Therefore, (\ref{dis1}) is
GAS, by Lemma \ref{impl}. The assertion in the theorem about
periodicity follows from Lemma \ref{fubl} and the formula for $U$,
so this proves our theorem.

\subsubsection{Proof of Theorem \ref{matrthm}}
\label{matrp}

Let $V_3 = V_1 + V_2$, and let $\alpha_1,\alpha_2\in {\mathcal
K}_\infty$  satisfy the UPPD requirements for $V_1$. In the rest of
the proof, all inequalities should be interpreted as holding
globally unless otherwise indicated. We also leave out the argument
$(x,k)$ of some of our $\mathcal{USB}$ functions when this would not
lead to confusion. It follows from Assumption \ref{as1} that we can
determine a positive definite function $\lambda$ such that
\begin{equation}
\label{lamjob} \Delta_k V_3(x,k)  \; \; \leq\; \;   - p(k+1) W(x) +
\phi_1(|x|)\phi_2(N_1(x,k))\; \; \le \; \;  - p(k+1)
\lambda(V_1(x,k)) + \phi_1(|x|)\phi_2(N_1(x,k))
\end{equation}
e.g. $\lambda(s)=\min\{W(x): x\in \mathbb{R}^n, \alpha_1(|x|)\le
s\le \alpha_2(|x|)\}$ (which is positive definite because $W\in
\mathcal{PD}$). By minorizing $\lambda$ as necessary as in
\cite{MM06}, we can assume it is $C^1$, nondecreasing on $[0,1/2]$
and nonincreasing on $[1/2,\infty)$. The proof of Lemma \ref{tl}
above with $\Theta(r):=\lambda(r/2)$ provides an increasing
continuous function $k_1:[0,\infty)\to[1,\infty)$ such that
$\Lambda_1(s):=k_1(s)\lambda(s)$ is of class ${\mathcal K}_\infty$.
Let $V_4 = k_1(V_1) V_3$. Then (\ref{lamjob}) gives
\begin{equation}
\label{mat7}\renewcommand{\arraystretch}{2}
\begin{array}{rcl}
\Delta_k V_{4}(x,k) & = & \left[k_1(V_1(F(x,k),k+1)) -
k_1(V_1(x,k))\right] V_{3}(F(x,k),k+1)
\\
& & + k_1(V_1(x,k))\left[V_{3}(F(x,k),k+1) - V_3(x,k)\right]
\\
& \leq & \left[k_1(V_1(F(x,k),k+1)) - k_1(V_1(x,k))\right]
V_{3}(F(x,k),k+1)
\\
& & - k_1(V_1(x,k)) p(k+1) \lambda(V_1(x,k))
 + k_1(V_1(x,k))\phi_1(|x|)\phi_2(N_1(x,k))\, .
\end{array}
\end{equation}
Since $F, V_3\in \mathcal{USB}$ and $V_1\in \mathcal{UPPD}$, we get
continuous increasing positive functions $\Gamma$ and $\Lambda_2$
such that
\begin{equation}
\label{mat8}
\begin{array}{rcl}
\Delta_k V_{4}(x,k) & \leq & [-\Delta_k V_{1}(x,k)] \Gamma(V_1(x,k))
- p(k+1) \Lambda_1(V_1(x,k))
 + \Lambda_2(V_1(x,k))\phi_2(N_1(x,k))
\end{array}
\end{equation}
e.g. by first finding an increasing positive function $\tilde
\alpha$ such that $|k'_1(r)|\le \tilde \alpha(r)$. Define $k_2\in
\mathcal{K}_\infty$ by $k_2(s)=s\Gamma(s)$. Since $k_2(s)/s$ is
increasing, we have $k_2(b)-k_2(a)\ge (b-a)\Gamma(b)$ when $b\ge
a\ge 0$ (since
 $(k_2(b)-k_2(a))(b-a)^{-1}\ge
k_2(b)b^{-1}$ when $b>a\ge 0$).   Hence by choosing
$a=V_1(F(x,k),k+1)$ and $b=V_1(x,k)$, we get
 $\Delta_k (k_2\circ V_1)(x,k)\le
\Delta_kV_1(x,k)\Gamma(V_1(x,k))$ everywhere.  Therefore, by adding
a $\mathcal{K}_\infty$ function to $k_2$ as necessary, we can assume
$V_5 := V_4 + k_2(V_1)\in \mathcal{UPPD}$ and satisfies
\begin{equation}
\label{mat9}
\begin{array}{rcl}
\Delta_k V_{5}(x,k) & \leq & - p(k+1) \Lambda_1(V_1(x,k)) +
\Lambda_2(V_1(x,k))\phi_2(N_1(x,k))\, .
\end{array}
\end{equation}
Arguing as in the proof of Theorem \ref{disthm} except with $\gamma$
replaced by $\Lambda_1$
 provides
$V_6\in\mathcal{UPPD}$  such that
\begin{equation}
\label{mat10}
\begin{array}{rcl}
\Delta_k V_{6}(x,k) \; \;  \leq \; \;  -
\displaystyle\frac{\delta}{4(l+1)}\Lambda_1(V_1(x,k)) +
\Lambda_2(V_1(x,k))\phi_2(N_1(x,k)). \end{array}\end{equation}

By arguing as in \cite[Section IV.A]{MM06}, we can select  $k_3\in
C^1\cap\mathcal{PD}$ so that $|k'_3(s)|\le 1$ for all $s\ge 0$ and
\begin{equation}\label{k3p}
k_3(r)\le
\phi^{-1}_2\left(\frac{\delta}{8(l+1)}\frac{\Lambda_1(r)}{1+\Lambda_2(r)}\right)\frac{1}{1+\Lambda_2(r)},\;
\; {\rm hence}\; \;
\phi_2\left(k_3(V_1)\Lambda_2(V_1)\right)\Lambda_2(V_1) \leq
\frac{\delta}{8(l+1)}\Lambda_1(V_1)\end{equation} everywhere. Choose
$\mu_F, \alpha_6\in \mathcal{K}_\infty$ such that $V_6(x,k)\le
\alpha_6(|x|)$ and $|F(x,k)|\le \mu_F(|x|)$ everywhere. Arguing as
above (with $\Gamma$ replaced by
$\alpha_6\circ\mu_F\circ\alpha^{-1}_1$) and recalling that
$\Delta_kV_1(x,k)\le 0$ provides $k_4\in \mathcal{K}_\infty$ such
that
\[
\Delta_k(k_4\circ V_1)(x,k) \; \le\;
\alpha_6\circ\mu_F\circ\alpha^{-1}_1\circ V_1(x,k)
\Delta_kV_1(x,k)\; \le\; \alpha_6(\mu_F(|x|))\Delta_kV_1(x,k)
\] everywhere.  Since
$|k_3'(s)|\le 1$ for all $s\ge 0$, we have \[[\Delta_k(k_3\circ
V_1)(x,k)]V_6(F(x,k), k+1)\; \le\;  [-\Delta_kV_1(x,k)]V_6(F(x,k),
k+1)\; \le\;  -\alpha_6(\mu_F(|x|))\Delta_kV_1(x,k)\] everywhere. It
follows that
 $V_7 := k_3(V_1) V_6 + k_4(V_1)\in \mathcal{UPPD}$
satisfies
\begin{equation}
\label{mat12} \begin{array}{rcl}\renewcommand{\arraystretch}{1.5}
\Delta_k V_{7}(x,k) &=& [\Delta_k(k_3\circ
V_1)(x,k)]V_6(F(x,k),k+1)+k_3(V_1(x,k))\Delta_kV_6(x,k)+\Delta_k(k_4\circ
V_1)(x,k)\\[.5em]
& \leq & - \frac{\delta}{4(l+1)}k_3(V_1)\Lambda_1(V_1) +
k_3(V_1)\Lambda_2(V_1)\phi_2(N_1(x,k))\, .
\end{array}
\end{equation}

Next note that for all functions $\mu \in {\cal K}_{\infty}$,
\begin{equation}\label{7decay}
\Delta_k V_{7} \; \;  \leq \; \;
-\frac{\delta}{4(l+1)}k_3(V_1)\Lambda_1(V_1) +
\mu\left(k_3(V_1)\Lambda_2(V_1)\right)k_3(V_1)\Lambda_2(V_1) +
\mu^{-1}(\phi_2(N_1(x,k)))\phi_2(N_1(x,k)).\end{equation} (The fact
that $ab\le \mu(a)a+\mu^{-1}(b)b$ for all $a,b\ge 0$ follows by
separately considering the cases where $\mu(a)\ge b$ and
$\mu^{-1}(b)\ge a$.) Choosing $\mu = \phi_2$ in (\ref{7decay}) gives
\begin{equation}
\label{mat14}
\begin{array}{rcl}
\Delta_k V_{7} & \leq & -
\frac{\delta}{4(l+1)}k_3(V_1)\Lambda_1(V_1) +
\phi_2\left(k_3(V_1)\Lambda_2(V_1)\right)k_3(V_1)\Lambda_2(V_1)
+ N_1(x,k) \phi_2(N_1)\\[.5em]& \leq & -
\frac{\delta}{8(l+1)}k_3(V_1)\Lambda_1(V_1) + N_1(x,k) \phi_2(N_1)\;
\; \; ({\rm by\ } (\ref{k3p})).
\end{array}
\end{equation}
 Therefore, since
$N_1\in \mathcal{USB}$ and $\phi_2\in \mathcal{K}_\infty$,  a
suitable function $\phi_3\in{\cal K}_{\infty}$ gives
\begin{equation}
\begin{array}{rcl}
\Delta_k V_{7} & \leq & -
\displaystyle\frac{\delta}{8(l+1)}k_3(V_1)\Lambda_1(V_1) + N_1(x,k)
\phi_3(V_1)\, .
\end{array}
\end{equation}
Arguing as in the construction of $k_2$ above (but with $\Gamma$
replaced by $\phi_3$) provides $k_5\in \mathcal{K}_\infty$ such that
$\Delta_k(k_5\circ V_1)(x,k)\le \Delta_kV_1(x,k)\phi_3(V_1(x,k))\le
-N_1(x,k)\phi_3(V_1(x,k))$.  Hence,
 $V_{8} := V_7 + k_5(V_1)\in \mathcal{UPPD}$ satisfies
\begin{equation}
\label{mat15}
\begin{array}{rcl}
\Delta_k V_{8}(x,k) & \leq & -
\displaystyle\frac{\delta}{8(l+1)}k_3(V_1(x,k))\Lambda_1(V_1(x,k))\;
\; \le \; \; -\alpha_3(|x|),
\end{array}
\end{equation}
where
\[\alpha_3(s)\; :=\; \frac{\delta}{8(l+1)}\min\{k_3(u)\Lambda_1(u):
\alpha_1(s)\le u\le\alpha_2(s)\}.\] Since $\alpha_3\in
\mathcal{PD}$, $V_8$ satisfies the requirements of the theorem. This
and Lemma \ref{impl} proves the theorem.

\subsection{Results for Hybrid Systems}
\label{mixes}

\subsubsection{Proof of Theorem \ref{hybthm}}

 For each
$k\in {\mathbb Z}_{\ge 0}$, let $V_{\rm cts}(x,t,k)$ denote the
continuous-time strictification of $V$ obtained in \cite{MM05} for
the nonstrictness parameter $q\in {\mathcal P}_{\rm cts}$.  Thus,
\[
V_{\rm
cts}(x,t,k):=\left[1+\frac{1}{\tau}\int_{t-\tau}^t\int_z^tq(\nu)\,
d\nu\, dz\right]V(x,t,k)\, .
\]
The results from \cite{MM05}   show that $\mathcal{D}V_{\rm
  cts}(x,t,k)\le -(\varepsilon/\tau) V_{\rm cts}(x,t,k)$ for all $x\in
  C$, $t\ge 0$, and $k\in \mathbb{Z}_{\ge 0}$.  We next rewrite
the first decay condition in (\ref{hyply2}) as
\begin{equation}\label{rw}
V(F(x,k),t,k+1) - V(x,t,k) \; \leq\;  -p(k+1) V(x,t,k),\; \; \;
\forall x\in D, t\in [0,\infty), k\in \mathbb{Z}_{\ge 0}
\end{equation}
where $k\mapsto p(k):=1-e^{-r(k)}$ is again of PE type. For each
$t\ge 0$, let $V_{\rm dis}(x,t,k)$ be the  strictification
\[
V_{\rm dis}(x,t,k):=\left[1+ \frac{1}{4(l+1)}\displaystyle\sum_{s =
k - l}^{k} \displaystyle\sum_{j = s}^{k}\, p(j)\right]V(x,t,k)
\]
of $V$ from Theorem \ref{disthm}.  The proof of  Theorem
\ref{disthm} shows we can take $\kappa(s)\equiv \gamma(s)\equiv s$,
and therefore also
\[ \Delta_kV_{\rm dis}(x,t,k)\; =\; V_{\rm
dis}(F(x,k),t,k+1) - V_{\rm dis}(x,t,k) \; \leq\;
-\frac{\delta}{4(l+1)}V_{\rm dis}(x,t,k)\; \; \; \forall t\in
[0,\infty), \, k\in \mathbb{Z}_{\ge 0}
\]
for all $x\in D$.  By enlarging $l$ from the PE assumption as
necessary, we can assume $\delta<l$.
 It follows that  the discrete decay condition
in (\ref{hyply2}) holds with $V$ replaced by $V_{\rm dis}$ and with
the constant
\[r(k)\equiv \ln\left(\frac{ 4(l+1)}{4(l+1)-\delta}\right)>0.\]
Since $\mathcal{D}V\le 0$ on $C$ and $\Delta_kV\le 0$ on $D$, we
have $\mathcal{D}V_{\rm dis}\le 0$ on $C$ and $\Delta_kV_{\rm
cts}\le 0$ on $D$. The uniform boundedness of  $S(k)$ and
$\int_{t-\tau}^t\int_z^tq(\nu)\, d\nu\, dz$ from Lemma \ref{fubl}
provides constants $r_c,r_d>0$ such that \[V_{\rm cts}(x,t,k)\le
r_cV_{\rm dis}(x,t,k)\le r_dV_{\rm cts}(x,t,k)\] everywhere. One
therefore easily checks that $V^\sharp(x,t,k):=V_{\rm
cts}(x,t,k)+V_{\rm dis}(x,t,k)$ as given by (\ref{vs}) is an
exponential decay Lyapunov function for the hybrid dynamic ${\cal
H}$.  The periodicity assertion follows as before from Lemma
\ref{fubl}, so the result follows from Lemma \ref{implh}.

\subsubsection{Proof of Theorem \ref{hmatrthm}}
For each $t\ge 0$ we apply the first part of the proof of Theorem
\ref{matrthm} to the functions $(x,k)\mapsto V_1(x,t,k)$ and
$(x,k)\mapsto V_2(x,t,k)$ to get $V_5$ that satisfies
\begin{equation} \label{mat9a}
\begin{array}{rcl}\! \! \! \! \Delta_k V_{5}(x,t,k) & \leq & - p(k+1)
\Lambda_1\left(V_1(x,t,k)\right) +
\Lambda_2\left(V_1(x,t,k)\right)\phi_2(N_1(x,t,k))\; \; \; \forall
x\in D, t\ge 0, k\in {\mathbb Z}_{\ge 0}.
\end{array}\end{equation}
This can be done with $\Lambda_1\in C^1$ and $\Lambda_2$ independent
of $t$.  For each $k\in {\mathbb Z}_{\ge 0}$, we next apply the
continuous time analog of the preceding argument (which is almost
exactly the same except with $\Delta_kV_i$ replaced by ${\mathcal
D}V_i$ for $i=1,2,\ldots, 5$, as discussed in the appendix below) to
get a continuous version $V^{\rm cts}_5$ of $V_5$ that satisfies
\begin{equation}
\label{mat9aa}
\begin{array}{rcl}
{\mathcal D}V^{\rm cts}_{5}(x,t,k) & \leq & - q(t)
\Lambda_1\left(V_1(x,t,k)\right) +
\Lambda_2\left(V_1(x,t,k)\right)\phi_2(N_1(x,t,k))\; \; \; \forall
x\in C, t\ge 0, k\in {\mathbb Z}_{\ge 0}.
\end{array}\end{equation}
In fact, by enlarging $k_2$  as necessary (e.g., by enlarging
$\Gamma$ in the discrete version of the proof), we can assume
$V^{\rm cts}_5$ and $V_5$ have the same formula. Applying the
strictification method from Theorem \ref{disthm}  to $V_5$ produces
\[
V^{\rm
dis}_6(x,t,k):=V_5(x,t,k)+\frac{1}{4(l+1)}S(k)\Lambda_1(V_5(x,t,k))
\]
that satisfies (\ref{mat10}) with $V_6$ replaced by $V^{\rm dis}_6$
and with   $V_1$ and $N_1$ now also depending on $t$.  Similarly, we
apply the continuous time strictification from \cite{MM05} (as in
the proof of Theorem \ref{hybthm}) to $V^{\rm cts}_5$ to get \[
V^{\rm cts}_6(x,t,k):= V^{\rm
cts}_5(x,t,k)+\frac{1}{\tau}\left[\int_{t-\tau}^t\int_z^tq(\nu)\,
d\nu\, dz\right]\Lambda_1(V^{\rm cts}_5(x,t,k)) \] that satisfies
${\mathcal D}V^{\rm cts}_6(x,t,k)\; \le\;
-\Lambda_1\left(V_1(x,t,k)\right)+\Lambda_2\left(V_1(x,t,k)\right)\phi_2(N_1(x,t,k))$
when $x\in C$, possibly by reducing $\Lambda_1$ and increasing
$\Lambda_2\in \mathcal{K}_\infty$ without relabeling. Setting
$V_6=V^{\rm cts}_6+V^{\rm dis}_6$, and assuming without loss of
generality that $1>\delta/\{4(l+1)\}$ (by enlarging $l$ without
relabeling as before), it follows from the fact that $V^{\rm cts}_5$
and $V_5$ have the same formula that we can enlarge $\Lambda_2$
sufficiently so that
\begin{eqnarray}
\label{mat10a} \Delta_k V_{6}(x,t,k) & \leq & -
\frac{\delta}{4(l+1)}\Lambda_1\left(V_1(x,t,k)\right) +
\Lambda_2\left(V_1(x,t,k)\right)\phi_2(N_1(x,t,k))\; \; \; \forall x\in D\\
{\mathcal D}V_{6}(x,t,k) & \leq & -
\frac{\delta}{4(l+1)}\Lambda_1\left(V_1(x,t,k)\right) +
\Lambda_2\left(V_1(x,t,k)\right)\phi_2(N_1(x,t,k))\; \; \; \forall
x\in C
\end{eqnarray}
hold for all $t\ge 0$ and $k\in {\mathbb Z}_{\ge 0}$.  (This can be
seen by bounding $|\Lambda'_1|$ on the relevant intervals and
recalling that $F,G\in \mathcal{USB}$ and $V_1\in \mathcal{UPPD}$.
In particular, to get (\ref{mat10a}), we write
$\Delta_k(\Lambda_1\circ V_5)(x,t,k)=\Lambda'_1(\eta
V_5(F(x,k),t,k+1)+(1-\eta)V_5(x,t,k))\Delta_kV_5(x,t,k)$ for
$\eta\in [0,1]$ depending on $x$, $t$,  and $k$ and use the fact
that $\Lambda'_1\ge 0$ everywhere.) We next follow the reminder of
the proof of Theorem \ref{matrthm} applied to $V_6$ for each $t\ge
0$ to get a function $V^{\rm dis}_8(x,t,k)$ satisfying the
conclusion of the proof when $x\in D$. We also apply the continuous
time analog of that part of the proof to $V_6$ for each $k$ (with
$\Delta_kV_i$ replaced by ${\mathcal D}V_i$ for all $i$ as before,
similarly to the argument done in the appendix below) to get $V^{\rm
cts}_8$ satisfying
\[{\mathcal D}V^{\rm cts}_8(x,t,k)\le -\tilde \alpha(|x|),\; \;
\forall x\in C, t\ge0, k\in {\mathbb Z}_{\ge 0}\] for a suitable
$\tilde \alpha\in \mathcal{PD}$. By enlarging $k_4, k_5\in {\mathcal
K}_\infty$ and reducing $k_3\in \mathcal{PD}$ in the continuous and
discrete versions of the proof, we can assume they are the same in
both versions, so $V^{\rm cts}_8$ and $V^{\rm dis}_8$ have the same
expression.  Hence,  we can satisfy the requirements of the theorem
with their common value. Combined with the result of Lemma
\ref{implh}, this proves the theorem.

\section{Examples}
\label{sec:examples}

One class of systems covered by our discrete time results is as
follows. Assume (\ref{dis1}) is GAS and that a strict Lyapunov
function $V$ for the system is available.  This provides
$\alpha_1,\alpha_2\in \mathcal{K}_\infty$ and $\alpha_3\in
\mathcal{PD}$ such that $\Delta_kV(x,k)\le -\alpha_3(|x|)$ and
$\alpha_1(|x|)\le V(x,k)\le \alpha_2(|x|)$ everywhere. Assume now
that the system is acted on by a PE term $p\in {\mathcal P}_{\rm
dis}$ that freezes the dynamics for certain times. The new system
becomes
\begin{equation}\label{fp}x_{k+1}=[1-p(k+1)]x_k+p(k+1)F(x_k,k).\end{equation}
Thus the new dynamic $F_p(x,k):=[1-p(k+1)]x+p(k+1)F(x,k)$ fixes the
state when $p(k+1)=0$.  By separately considering the cases
$p(k+1)=0$ and $p(k+1)=1$, one checks that if $p(k)\in \{0,1\}$ for
all $k$, then $V(F_p(x,k),k+1)-V(x,k)\le -p(k+1)\alpha_3(|x|) \le
-p(k+1)\Theta(V(x,k))$ everywhere, where
$\Theta(s)=\min\{\alpha_3(p): \alpha^{-1}_2(s)\le p\le
\alpha^{-1}_1(s)\}$. Since $\Theta\in \mathcal{PD}$, $V$ satisfies
the PE decay condition from Theorem \ref{disthm} for the new dynamic
$F_p$. More generally, assume $p(k)\in [0,1]$ for all $k$. Assume
also that $V(x,k)$ is a Lyapunov function for (\ref{dis1}) that is
independent of $k$ and convex in $x$.  Choose $\alpha_3\in
\mathcal{PD}$ such that $V(F(x,k))-V(x)\le -\alpha_3(|x|)$
everywhere. Then
\begin{eqnarray}
V(F_p(x,k))-V(x)&\le & [1-p(k+1)]V(x)+p(k+1)V(F(x,k))-V(x)\nonumber\\
&\le &-p(k+1)V(x)+p(k+1)[V(x)-\alpha_3(|x|)]\; \; = \; \;
-p(k+1)\alpha_3(|x|)\nonumber
\end{eqnarray}
everywhere, so $F_p$ again satisfies our PE assumptions.

A general class of {\em hybrid} systems covered by our
strictification results is as follows.  Assume the continuous time
system (\ref{con1}) admits $q\in \mathcal{P}_{\rm cts}$, $\gamma\in
\mathcal{K}_\infty$, $V\in C^1$, and  $\alpha_1,\alpha_2\in
{\mathcal K}_\infty$ satisfying $\mathcal{D}V(x,t)\le
-q(t)\gamma(V(x,t))$ and $\alpha_1(|x|)\le V(x,t)\le \alpha_2(|x|)$
for all $x\in \mathbb{R}^n$ and $t\ge 0$ (i.e., (\ref{con1}) admits
a nonstrict Lyapunov function in the sense of \cite{MM05}).
 \footnote{A concrete example where this occurs and where it is easy
to find $V$ is where $q\in \mathcal{P}_{\rm cts}$ (e.g.
$q(t)=\sin^2(t)$) and $\dot x=h(x,t)$ is GAS (e.g. $\dot x=-x$) and
we take the dynamic $G(x,t)=q(t)h(x,t)$ and a Lyapunov function
$V(x,t)$ for $\dot x=h(x,t)$.}
 Given subsets $C, D\subseteq \mathbb{R}^n$
and $p\in \mathcal{P}_{\rm dis}$ taking all its values in $\{0,1\}$,
we determine conditions on $F\in \mathcal{USB}$ guaranteeing that we
can construct a Lyapunov function for
\begin{equation}
\label{hybsyp} {\mathcal H}_p:= \left\{
\begin{array}{lcll}
\dot x&=&G(x,t)\; \; , & x\in C\\
x_{k+1}&=&F_p(x_k,k)\; \; ,&  x_k\in D
\end{array},
\right.\end{equation} where $F_p$ is as defined  above. (The
construction we are about to give also works if instead of assuming
$p(k)\in \{0,1\}$ for all $k\in \mathbb{Z}_{\ge 0}$, we assume (i)
$x\mapsto V(x,t)$ is convex for each $t\in [0,\infty)$ and (ii)
$p(k)\in [0,1]$ for all $k\in \mathbb{Z}_{\ge 0}$.  This situation
arises if $\dot x=G(x,t):=A(t)x$ is GAS and $A(t)$ is continuous and
bounded since then we can take $V(x,t):=x^\top P(t)x$ for a suitable
bounded everywhere positive definite matrix $P(t)$ \cite[Section
4.6]{K02}.)
 To this end,
 first notice that by reducing
$\gamma\in \mathcal{K}_\infty$ as necessary, we can assume
$\gamma\in C^1$ and $\gamma(s)\le \alpha_1(\alpha^{-1}_2(s))/2$ for
all $s\ge 0$. Let $F$ satisfy $|F(x,k)|\le
\alpha^{-1}_2(\alpha_1(|x|)/2)$ for all $x\in D$ and $k\in
\mathbb{Z}_{\ge 0}$.  (This reduces to a linear growth condition
when $V(x,t)=x^\top P(t) x$ and $P$ has bounded positive
eigenvalues.) By separately considering the cases $p(k+1)=0$ and
$p(k+1)=1$, it follows that
\begin{eqnarray*}
V(F_p(x,k),t)-V(x,t)&\le&
p(k+1)\alpha_2(|F(x,k)|)-p(k+1)\alpha_1(|x|)\\&\le&
-\frac{1}{2}p(k+1)\alpha_1(|x|)\; \; \le \; \;
-p(k+1)\gamma(\alpha_2(|x|))\; \; \le\; \;
-p(k+1)\gamma(V(x,t))\end{eqnarray*} for all $x\in D$, $t\ge 0$, and
$k\in \mathbb{Z}_{\ge 0}$.  A slight variant of the proof of Theorem
\ref{hybthm}  therefore provides an explicit globally smooth strict
Lyapunov function for ${\mathcal H}_p$ having the form
\[
V^\sharp(x,t,k):=2V(x,t)+\frac{1}{\tau}\left[\int_{t-\tau}^t\int_s^tq(r)\,
dr\, ds\right]\gamma(V(x,t))+\left[\frac{1}{4(l+1)}\sum_{s = k -
l}^{k} \sum_{j = s}^{k}\, p(j)\right]\gamma(V(x,t))
\]
for  $l$ and $\tau$ as in the requirements $p\in \mathcal{P}_{\rm
dis}$ and $q\in \mathcal{P}_{\rm cts}$ so ${\mathcal H}_p$ is  GAS,
as claimed.

\section{Conclusions}
\label{sec:6}

We provided new methods for constructing closed form strict Lyapunov
functions for hybrid systems that admit appropriate nonstrict
Lyapunov functions.  Our results cover cases where the given
nonstrict Lyapunov functions satisfy  a decay condition involving
persistency of excitation  parameters or hybrid versions of the
conditions of Matrosov's Theorem.  Due to the ubiquity of Lyapunov
functions in engineering applications, we expect that our results
will be useful in a wide range of settings in which explicit
Lyapunov functions are needed such as Lyapunov-based controller
design and robustness analysis.  We conjecture that our results can
be extended to hybrid control systems with outputs. This would
extend  \cite{SW99, SW01} and the input-to-output stability Lyapunov
function constructions from \cite{MM05b} to hybrid systems and also
provide more explicit constructions that would be suited for
applications.

 \renewcommand{\theequation}{A.\arabic{equation}}
\renewcommand{\thetheorem}{A.\arabic{theorem}}
  \setcounter{equation}{0}
\setcounter{theorem}{0}

  \section*{Appendix}
Our proof of Theorem \ref{hmatrthm} was based on a continuous time
version of Theorem \ref{matrthm}.  We next give a precise statement
and proof of  this continuous time result, which is applied in the
proof of Theorem \ref{hmatrthm} to $(x,t)\mapsto V_i(x,t,k)$ for
each $k$ and $i=1,2$.
  We assume the following version of the Matrosov conditions:
\begin{assumption}\label{as1aa}  There exist
$V_1:{\mathbb R}^n\times [0,\infty)\to [0,\infty)$ of class
$\mathcal{UPPD}$ and $V_2:{\mathbb R}^n\times [0,\infty)\to {\mathbb
R}$ of class $\mathcal{USB}$ that are $C^1$, $\phi_2\in {\mathcal
K}_\infty$, nonnegative functions $N_1, N_2\in \mathcal{USB}$, a
function $\chi:{\mathbb R}^n\times [0,\infty)\times [0,\infty)\to
{\mathbb R}$, a positive increasing function $\phi_1$, $W\in
\mathcal{PD}$, and $q\in {\cal P}_{\rm cts}$ such that
\[\begin{array}{l}
 \mathcal{D} V_{1}(x,t) \leq - N_1(x,t),\; \; \;\; \mathcal{D}
V_{2}(x,t) \leq - N_2(x,t) + \chi(x,N_1(x,t),t),\\[.5em]
 |\chi(x,N_1(x,t),t)| \leq
\phi_1(|x|)\phi_2(N_1(x,t)),\; \; {\rm and}\; \; N_1(x,t) + N_2(x,t)
\geq q(t) W(x)
\end{array}\]
hold  for all $x\in {\mathbb R}^n$ and $t\in [0,\infty)$.
\end{assumption}
Notice that $V_2$ can take both positive and negative values. We
show:

\begin{theorem}
\label{cmatrthm} If (\ref{con1}) satisfies Assumption \ref{as1aa},
then one can construct an explicit strict Lyapunov function for
(\ref{con1}).  In particular, (\ref{con1}) is GAS.\end{theorem}

 To prove this theorem, we
indicate the changes needed in the  proof of Theorem \ref{matrthm}.
 We define $V_3$ and $\lambda$
as in Section \ref{matrp} which therefore satisfy
\[
\mathcal{D} V_3(x,t)  \; \; \leq\; \;   - q(t) W(x) +
\phi_1(|x|)\phi_2(N_1(x,t))\; \; \le \; \;  -q(t) \lambda(V_1(x,t))
+ \phi_1(|x|)\phi_2(N_1(x,t))
\]
everywhere.  We also define $k_1$, $\Lambda_1$, and
$V_4:=k_1(V_1)V_3$ as before and as before also determine positive
increasing functions $\Gamma$  and $\Lambda_2$ such that
\begin{equation}
\label{mat8a}
\begin{array}{rcl}
\mathcal{D} V_{4}(x,t) & \leq & [-\mathcal{D} V_{1}(x,t)]
\Gamma(V_1(x,t)) - q(t) \Lambda_1(V_1(x,t))
 + \Lambda_2(V_1(x,t))\phi_2(N_1(x,t)).
\end{array}
\end{equation}
Choosing $k_2\in \mathcal{K}_\infty$ such that $k'_2\ge \Gamma$
everywhere gives $\mathcal{D}(k_2\circ
V_1)=k'_2(V_1)\mathcal{D}V_1\le \Gamma(V_1)\mathcal{D}V_1$, since
$\mathcal{D}V_1\le 0$ everywhere.  Enlarging $k_2\in
\mathcal{K}_\infty$ as necessary, it follows that
$V_5:=V_4+k_2(V_1)\in \mathcal{UPPD}$ satisfies
\begin{equation}
\label{mat9b}
\begin{array}{rcl}
\mathcal{D} V_{5}(x,t) & \leq & - q(t) \Lambda_1(V_1(x,t)) +
\Lambda_2(V_1(x,t))\phi_2(N_1(x,t))\, .
\end{array}
\end{equation}
Applying the continuous time strictification method of \cite{MM05}
and enlarging $\Lambda_2$ and reducing $\Lambda_1$ as necessary
without relabeling provides $\gamma\in \mathcal{K}_\infty$ and
$\tau>0$ such that
\begin{equation}\label{new6}
V_6(x,t):=V_5(x,t)+\left[\int_{t-\tau}^t\int_s^tq(r)\, dr\,
ds\right]\, \gamma(V_5(x,t))
\end{equation}
satisfies $\mathcal{D} V_{6}(x,t)  \leq  -  \Lambda_1(V_1(x,t)) +
\Lambda_2(V_1(x,t))\phi_2(N_1(x,t))$.  This uses the global
boundedness of the double integral in (\ref{new6}) from Lemma
\ref{fubl}.  Arguing as in the proof of Theorem \ref{matrthm} gives
 $k_3\in
\mathcal{PD}\cap C^1$ such that
\[
k_3(r)\le
\phi^{-1}_2\left(\frac{\Lambda_1(r)}{1+\Lambda_2(r)}\right)\frac{1}{1+\Lambda_2(r)},\;
\; {\rm hence}\; \;
\phi_2\left(k_3(V_1)\Lambda_2(V_1)\right)\Lambda_2(V_1) \leq
\Lambda_1(V_1)
\]
everywhere.  Choose $k_4\in \mathcal{K}_\infty\cap C^1$ such that
$k'_4(s)\ge |k'_3(s)|(\alpha_6\circ \alpha^{-1}_1)(s)$ everywhere,
where $\alpha_1$ is as in the UPPD requirement on $V_1$, and
$\alpha_6(|x|)\ge V_6(x,t)$ for all $x\in \mathbb{R}^n$ and $t\ge
0$.  Then $k'_4(V_1)\ge |k'_3(V_1)|V_6$ everywhere, so
$V_7:=k_3(V_1)V_6+k_4(V_1)\in \mathcal{UPPD}$ everywhere satisfies
\[
\begin{array}{rcl}
\mathcal{D}V_7&\le & -k_3(V_1)\Lambda_1(V_1)\; +\;
k_3(V_1)\Lambda_2(V_1)\phi_2(N_1)
\; -\; V_6|k'_3(V_1)|\mathcal{D}V_1\; +\; k'_4(V_1)\mathcal{D}V_1\\
&\le & -k_3(V_1)\Lambda_1(V_1)\; +\;
k_3(V_1)\Lambda_2(V_1)\phi_2(N_1).
\end{array}
\]
The rest of the argument is similar to the corresponding part of the
proof of Theorem \ref{matrthm} with $\Delta_k$ replaced by $\mathcal
D$ and $\frac{\delta}{8(l+1)}$ replaced by $1$.


\begin{thebibliography}{99}

\bibitem{AlS99}
 Albertini, F., and E.D. Sontag, ``Continuous control-Lyapunov functions
 for asymptotically controllable time-varying systems,'' {\em International Journal of
  Control},  {\bf 72}(1999), pp. 1630-1641.
\bibitem{A99}
Angeli, D.,  ``Input-to-State Stability of PD-controlled robotic
systems,'' {\em Automatica}, {\bf 35}(1999), pp. 1285-1290.

\bibitem{AS99} Angeli, D., and E.D. Sontag, ``Forward
completeness, unboundedness observability, and their Lyapunov
characterizations,'' {\em Systems and Control Letters},  {\bf
38}(1999), pp. 209-217.


\bibitem{ASW00} Angeli, D., E.D. Sontag, and Y. Wang, ``A
characterization of integral input to state stability,'' {\em IEEE
Transactions on Automatic Control},  {\bf 45}(2000), pp. 1082-1097.

\bibitem{BR01}
Bacciotti, A., and L. Rosier, {\em Liapunov Functions and Stability
in Control Theory},  Springer, London, 2001.

\bibitem{CTG05} Cai, C., A. Teel, and R. Goebel, ``Converse
Lyapunov theorems and robust asymptotic stability for hybrid
systems,'' in {\em Proceedings of the 24th American Control
Conference (Portand, OR, June 2005)}, pp. 12-17,
\htmladdnormallink{http://www.ccec.ece.ucsb.edu/$\sim$cai/}
{http://www.ccec.ece.ucsb.edu/~cai/}.


\bibitem{C04} Collins, P., ``A trajectory-space approach to hybrid
systems,'' in {\em Proceedings of the International Symposium on the
Mathematical Theory of Networks and Systems (Katholiek Univ. Leuven,
Belgium, August 2004)}, Paper \#250,
\htmladdnormallink{http://homepages.cwi.nl/$\sim$collins/}{http://homepages.cwi.nl/~collins/}.

\bibitem{C92} Coron, J.-M., ``Global asymptotic stabilization for
controllable systems without drift,'' {\em Mathematics of Control,
Signals \& Systems}, {\bf 5}(1992), pp. 295-312.

\bibitem{FP00} Faubourg, L., and J.-B. Pomet, ``Control Lyapunov
functions for homogeneous ``Jurdjevic-Quinn'' systems,'' {\em ESAIM:
Control, Optimisation and Calculus of Variations}, {\bf 5}(2000),
pp. 293-311.


\bibitem{K02}
Khalil, H., {\it Nonlinear Systems, Third Edition}, Prentice Hall,
Englewood Cliffs, NJ,  2002.

\bibitem{KSW01} Krichman, M., E.D. Sontag, and Y. Wang,
``Input-output-to-state stability,'' {\em SIAM Journal on  Control
and Optimization}, {\bf 39}(2001), pp. 1874-1928.



\bibitem{MM05b} Malisoff, M., and F. Mazenc, ``Further
constructions of strict Lyapunov functions for time-varying
systems,'' in {\em Proceedings of the American Control Conference
(Portland, OR, June 2005)}, pp. 1889-1894.


\bibitem{MM05}
Malisoff, M., and F. Mazenc, ``Further remarks on strict
input-to-state stable Lyapunov functions for time-varying systems,''
{\em Automatica}, {\bf 41}(2005), pp. 1973-1978.



\bibitem{MGSW05a} Mancilla-Aguillar, J., R. Garcia, E.D. Sontag,
and Y. Wang, ``On the representation of switched systems with inputs
by perturbed control systems,''  {\em Nonlinear Analysis:  Theory,
Methods \& Applications}, {\bf 60}(2005), pp. 1111-1150.

\bibitem{MGSW05b}
Mancilla-Aguillar, J., R. Garcia, E.D. Sontag, and Y. Wang,
``Uniform stability properties of switched systems with switchings
governed by digraphs,'' {\em Nonlinear Analysis:  Theory, Methods \&
Applications}, {\bf 63}(2005), pp. 472-490.


\bibitem{M03}
Mazenc, F., ``Strict Lyapunov functions for time-varying systems,''
{\em Automatica},  {\bf 39}(2003), pp. 349-353.

\bibitem{MM06} Mazenc, F., and M. Malisoff, ``Further
constructions of control-Lyapunov functions and stabilizing
feedbacks for systems satisfying the Jurdjevic-Quinn conditions,''
{\em IEEE Transactions on Automatic Control}, {\bf 51}(2006),  pp.
360-365.







\bibitem{MN06} Mazenc, F.,  and D. Nesi\'c, ``Lyapunov functions
for time varying systems satisfying generalized conditions
 of Matrosov theorem,'' in {\em Proceedings of the
44th IEEE Conference on Decision
 \& Control (CDC) \& European Control Conference ECC 05 (Seville, Spain, December
 2005)},   pp. 2432-2437.

\bibitem{MSP99} Morin,  P.,  C. Samson, and J.-B. Pomet, ``Design
of homogeneous time-varying stabilizing control laws for driftless
systems via oscillatory approximation of Lie brackets in closed
loop,'' {\em SIAM Journal  on Control and Optimization}, {\bf
38}(1999), pp. 22-49.


 \bibitem{NL04}
  Nesi\'c, D.,  and A. Loria, ``On uniform asymptotic stability of time-varying
  parameterized discrete-time cascades,'' {\em IEEE Transactions on Automatic
  Control}, {\bf 49}(2004), pp. 875-887.

  \bibitem{NTS99}
  Nesi\'c, D., A. Teel, and E.D. Sontag, ``Formulas relating
  $\mathcal{KL}$ stability estimates of discrete-time and
  sampled-time nonlinear systems,'' {\em Systems and Control
  Letters}, {\bf 38}(1999), pp. 49-60.

  \bibitem{Sa90}
  Samson, C.,
``Velocity and torque feedback control of a nonholonomic cart,'' in
   {\em Advanced Robot Control (Grenoble, 1990)},
 Lecture Notes in Control and Information  Sciences Vol. 162, Springer, Berlin,
 1991, pp. 125-151.

\bibitem{S89}
Sontag, E.D.,   ``Smooth stabilization implies coprime
factorization,'' {\em IEEE Transactions on   Automatic   Control},
{\bf 34}(1989), pp. 435--443.

\bibitem{So90}
Sontag, E.D., ``Feedback stabilization of nonlinear systems,'' in
{\em Robust Control of Linear Systems and Nonlinear Control} (M.A.
Kaashoek, Ed.), Birkh\"auser, Basel, 1990, pp. 61-81.

\bibitem{SW99} Sontag, E.D., and Y. Wang, ``Notions of input to
output stability,'' {\em Systems and Control Letters}, {\bf
38}(1999), pp. 235--248.

\bibitem{SW01} Sontag, E.D., and Y. Wang, ``Lyapunov
characterizations of input to output stability,'' {\em SIAM Journal
on  Control and Optimization}, {\bf 39}(2001), pp. 226-249.

\bibitem{VS00}
Van der Schaft, A., and H. Schumacher, {\em An Introduction to
Hybrid Systems}, Lecture Notes in Control and Information Sciences
Vol. 251, Springer-Verlag London, Ltd., London, 2000.

\end{thebibliography}
\end{document}